\documentclass{amsart}
\usepackage[utf8]{inputenc}
\usepackage{yufei}

\bibliography{citation.bib}

\title{Common kings of a chain of cycles in a strong tournament}

\author{Logan Post}
\address{Haverford College, Haverford, PA, 19041}
\email{lpost@haverford.edu}

\author{Zeyu Zheng}
\address{School of Mathematical Sciences, Fudan University, Shanghai, China 200433}
\email{zeyuzheng19@fudan.edu.cn}

\begin{document}
\begin{abstract}
    It is known that every strong tournament has directed cycles of any length, and thereby strong subtournaments of any size.
    In this note, we prove that they also can share a common vertex which is a king of all of them. This common vertex can be any king in the whole tournament. Further, the Hamiltonian cycles in them can be recursively constructed by inserting an additional vertex to one directed edge.
\end{abstract}
\maketitle

\section{Introduction}
A {\bf tournament} is a directed graph, obtained by assigning an orientation to each edge of a complete graph. A {\bf king} is a vertex $k$ such that each other vertex is reachable from $k$ by a directed path of length at most $2$ \cite{landau1953dominance}. A \textbf{strong} tournament has the property that every vertex is reachable from every other by a directed path \cite{moon1968topic}. A {\bf chain of cycles} is a sequence $C_1,C_2,\ldots,C_n$ such that for  $i\in[1,n-1]$, the cycle $C_{i+1}$ is obtained from $C_i$ by replacing some edge $(x,y)$ of  $C_i$ by two edges $(x,z)$ and $(z,y)$ where $z$ is not on $C_i$.

%\begin{definition}
%\label{def:bv}
%In a tournament $T_n$, we call a vertex $v$ a {\it Bottom vertex} if there exists a directed path from any vertex in $V(T_n)\setminus\{v\}$ to $v$.
%\end{definition}
%\begin{lemma}
%\label{lem:1}
%There exists a Hamiltonian path ending at any {\it Bottom vertex}.
%\end{lemma}
It is known that every tournament has a king \cite{landau1953dominance}, and it is known that every strong tournament $T_n$ has directed cycles of length $3,...,n$ \cite{moon1968topic}. The proof can be modified so that all cycles contain a common vertex which can be any king in the whole tournament, but not necessarily of the subtournaments determined by the vertices of the cycles.

In this note, we show that:
\begin{restatable}{theorem}{thmmain}
\label{thm:main}
 Given a strong tournament $T_n$ on $n$ vertices and any king $k$ of $T_n$, there exist a chain of cycles $C_3,C_4,\ldots,C_n$ where $C_i$ is of length $i$ such that $k$ is a king in all subtournaments induced by the vertex sets of the cycles in the chain.
\end{restatable}

%Equivalently, given a strong tournament $T_n$ and a king $k$, we can enumerate the vertices $k,v_1,...,v_{n-1}$ so for any integer $i\in[2,n-1]$, $\{k,v_1,...,v_i\}$ determines a strong subtournament where $k$ is the king.

\section{Proof of \texorpdfstring{\cref{thm:main}}{Theorem 1}}
\label{sec:proof}

Let $k$ be a king of $T_n$. We use $A$ to denote $N^+(k)$ and use $B$ to denote $N^-(k)$. A result in \cite{reid1978tournaments} shows that we can partition $A$ into $A_1,A_2,\ldots,A_r$ such that $A_i$ determines a strong subtournament for all integer $i\in[1,r]$, and for all $i<j$, all directed edges between $A_i$ and $A_j$ terminates in $A_j$. %We refer to $A_r$ as the root component in the decomposition.

\begin{claim} 
\label{edgegoesdown}
There exists $a^*\in A_r$, $b^*\in B$ such that $(a^*,b^*)$ is an edge in $T_n$.
\end{claim}

For a chosen $a\in A_r$, as $T_n$ is strong, there exists a directed path from $a$ to $k$. Choose the shortest such path. This path will look like $a\rightarrow \cdots a^*\rightarrow b^*\rightarrow k$ for some $b^*\in B$. $a^*\in A_r$ because $A_r$ is fully beaten by $A_i$ for all $i<r$, so this path resides in $A_r$ and no other elements of $B$ are included because this would imply a shorter path. This satisfies Claim \ref{edgegoesdown}. 

\begin{claim}
\label{hpthroughA}
There is a Hamiltonian path through $A$ ending at $a^*$.
\end{claim} 
We see that the subtournament on $A_1\cup A_2\cup...\cup A_{r-1}$ has Hamiltonian path $P_1$. Also, $A_r$ is strong, so it has a Hamiltonian cycle, in which there is a Hamiltonian path $P_2$ ending at $a^*$. Concatenating $P_1$ and $P_2$, this satisfies Claim \ref{hpthroughA}. We denote this path by $P=a_1\rightarrow a_2\rightarrow\cdots\rightarrow a_d=a^*$ where $d=\deg^+k$.
\iffalse
\begin{figure}[H]
\centering
\includegraphics[scale=.5]{Conj_5.png}
\caption{$T_n$, showing some known relationships within $A$ and $B$.}
\end{figure}
\fi

Using Claim \ref{edgegoesdown} and Claim \ref{hpthroughA}, we can construct length-$i$ directed cycles $C_i$ of which $k$ is the king for $i\in [3,n]$.

For any $i\in [1,d]$, we see that $C_{i+2}=k\rightarrow a_{d-i+1}\rightarrow \cdots\rightarrow a_d\rightarrow b\rightarrow k$ is a directed cycle of length $i+2$ and $k$ is the king of the subtournament on $V(C_{i+2})$. This gives us cycles of length from $3$ to $d+2$.

Given a directed cycle $C_{i-1}$ of length $i-1$ where $i\in[d+3, n]$ and $(A\cup\{k\})\subset V(C_{i-1})$. We now show that we can always find a cycle of length $i$ and its vertex set still contains $A\cup\{k\}$. 

Since $T_n$ is strong we know that there exists an edge $(w,z)\in T_n$ so that $w\in C_{i-1}$ and $z\in V(T_n)\setminus V(C_{i-1})$. Note that $(V(T_n)\setminus V(C_{i-1}))\subset B$, so $(z,k)$ is an edge in $T_n$. Therefore, there exists an edge $(x,y)$ on the cycle $C_{i-1}$ such that both $(x,z)$ and $(z,y)$ are edges in $T_n$. We can insert $z$ in the cycle $C_{i-1}$ between $x$ and $y$ to get a directed cycle $C_i$ of length $i$, and $(A\cup\{k\})\subseteq V(C_{i})$. This recursively gives us cycles of length from $d+3$ to $n$.

Now we have created a collection of cycles $C_3,C_4,\ldots ,C_n$ in $T_n$ of which $k$ is the king of subtournament induced by vertices of each of them. Note as an added fact that $C_{i+1}$ can be obtained from $C_i$ by replacing some edge $(x,y)$ of $C_i$ by two edges $(x,z)$ and $(z,y)$ where $z\in V(T_n)\setminus V(C_i)$ for all integer $i\in[3,n-1]$, which finishes the proof.

\section*{Acknowledgements}
The authors thank Professor Andr\'as Gy\'arf\'as for proposing the research problem we solve here in his course Advanced Combinatorics and for his valuable suggestions. The work have been conducted under the auspices of the Budapest Semesters in Mathematics Program during the Spring of 2022.

Yaobin Chen, Professor Jie Ma, Professor Hehui Wu have provided different proofs to the main theorem. The second author would like to thank them for fruitful early discussions on this problem.

\printbibliography

\end{document}